\documentclass[a4paper,11pt]{amsart}
\usepackage{tipa}
\usepackage{txfonts}
\usepackage{amssymb}
\usepackage{amsmath}
\usepackage{mathrsfs}
\usepackage{amsmath,amssymb,amsthm,latexsym,amscd,mathrsfs}
\usepackage{indentfirst}
\usepackage{graphicx}
\usepackage{booktabs}
\usepackage{array}
\usepackage{chemarrow}
\newcommand{\ROM}[1]{\mathrm{\uppercase\expandafter{\romannumeral#1}}}
\theoremstyle{definition}

\newtheorem{thm}{Theorem}[section]
\newtheorem{lem}{Lemma}[section]

\newtheorem{cor}{Corollary}[section]
\newtheorem{examp}{Example}[section]
\newtheorem{rem}{Remark}[section]
\newtheorem{prop}{Proposition}[section]

\newtheorem{prob}[thm]{Problem}

\newtheorem{ack}{Acknowledgements}   
\makeatletter


\title[Isoparametric functions on exotic spheres]{\textbf{Isoparametric functions on exotic spheres}}
\author[C. Qian]{Chao Qian}\address{School of Mathematical Sciences, University of Chinese Academy of Sciences , Beijing 100049, P.R.China; and
School of Mathematical Sciences, Laboratory of Mathematics and Complex Systems, Beijing Normal
University, Beijing 100875, P.R.China}\email{qianchao@ucas.ac.cn}
\thanks{The project is partially supported by the NSFC ( No.11071018
and No.11401560 )}
\author[Z. Z. Tang]{Zizhou Tang}\address{School of Mathematical Sciences, Laboratory of Mathematics and Complex Systems, Beijing Normal
University, Beijing 100875, P.R.China}\email{zztang@bnu.edu.cn}
\thanks{The second author is the corresponding author.}

\date{}
\begin{document}

\maketitle
\begin{abstract}
This paper extends widely the work in \cite{GT13}. Existence and
non-existence results of isoparametric functions on exotic spheres
and Eells-Kuiper projective planes are established. In particular,
every homotopy $n$-sphere ($n>4$) carries an isoparametric function
(with certain metric) with 2 points as the focal set, in strong
contrast to the classification of cohomogeneity one actions on
homotopy spheres \cite{St96} ( only exotic Kervaire spheres admit
cohomogeneity one actions besides the standard spheres ). As
an application, we improve a beautiful result of B\'{e}rard-Bergery
\cite{BB77} ( see also pp.234-235 of \cite{Be78} ).
\end{abstract}


\section{Introduction}

Let $N$ be a connected complete Riemannian manifold. A non-constant
smooth function $f$ on $N$ is called \emph{transnormal}, if there
exists a smooth function $b:\mathbb{R}\rightarrow\mathbb{R}$ such
that $|\nabla f|^2=b(f),$ where $\nabla f$ is the gradient of $f$.
If in addition, there exists a continuous function
$a:\mathbb{R}\rightarrow\mathbb{R}$ so that $\triangle f=a(f),$
where $\triangle f$ is the Laplacian of $f$, then $f$ is called
\emph{isoparametric}. Each regular level hypersurface is called an
\emph{isoparametric hypersurface} and the singular level set is
called the\emph{ focal set}. The two equations of the function $f$
mean that the regular level hypersurfaces of $f$ are parallel and have constant mean
curvatures, which may be regarded as a geometric generalization of
cohomogeneity one actions in the theory of transformation groups (
ref. \cite{GT13} ).

Owing to E. Cartan and H. F. M\"{u}nzner \cite{Mu80}, the
classification of isoparametric hypersurfaces in a unit sphere has
been one of the most challenging problems in submanifold geometry.
Up to now, the classification has almost been completed. We refer to
\cite{CR85} and \cite{Th00} for the development of this subject. For the most recent
progress and applications, see for example \cite{CCJ07}, \cite
{Miy13}, \cite{Miy14}, \cite{Ch13} and \cite{TY13}.

In general Riemannian manifolds, a series of beautiful results,
similar to the case in a unit sphere had been proved or claimed by
Wang in \cite{Wa87}. Based on Wang's work, Ge and Tang \cite{GT13}
made new contributions to the subject of isoparametric functions on
general Riemannian manifolds, especially on exotic spheres. In a
word, the results of \cite{Wa87} and \cite{GT13} show that the
existence of an isoparametric function ( even a transnormal function
) on a Riemannian manifold $N$ may restrict strongly the geometry
and topology of $N$. However, as posed in \cite{GT13}, it was not
known whether any exotic sphere with dimension more than $4$ admits
a metric and an isoparametric function with $2$ points as the focal
set ( one of the main results in \cite{GT13} asserts that an exotic
$4$-sphere ( if exists ) admits no properly transnormal functions ).
Moreover, a transnormal function $f$ on the Gromoll-Meyer $7$-sphere
with $2$ points as the focal set was constructed in \cite{GT13}, though
$f$ is not isoparametric.

One of our main aims is to develop a general way to construct
metrics and isoparametric functions on a given manifold. According
to \cite{Wa87}, one can see that a transnormal function on a
complete Riemannian manifold is necessarily a Morse-Bott function.
As is well known that a Morse-Bott function is a generalization of a
Morse function, and it admits critical submanifolds satisfying
certain non-degenerate condition. As one of our prime results we
show:

\begin{thm}\label{Main}
( A fundamental construction ) Let $N$ be a closed connected smooth
manifold and $f$ a Morse-Bott function on $N$ with the critical set
$C(f)=M_+\sqcup M_-$, where $M_+$ and $M_-$ are both closed
connected submanifolds of codimensions more than 1. Then there
exists a metric $g$ on $N$ so that $f$ is an isoparametric function.
Moreover, the metric can be chosen so that the critical submanifolds
${M}_+$ and ${M}_-$ are both totally geodesic.
\end{thm}
Recall that a closed smooth $n$-manifold is called a \emph{homotopy
$n$-sphere}, if it has the homotopy type of the unit sphere $S^n$.
By an \emph{exotic sphere} we mean a homotopy sphere not
diffeomorphic to $S^n$. By using a theorem of S. Smale, we have the
following result as a consequence.

\begin{cor}\label{homosphere}
Every homotopy $n$-sphere with $n>4$ admits a metric and an
isoparametric function with 2 points as the focal set.
\end{cor}

\begin{rem}
Corollary \ref{homosphere} answers partially Problem 4.1 in
\cite{GT13}.
\end{rem}

Moreover, we can also construct metrics and isoparametric functions
on homotopy spheres and the Eells-Kuiper projective planes so that
at least one component of the critical set is not a single point.

Next, we intend to ask if our proof of Theorem \ref{Main} can be
modified so that each regular level hypersurface has more geometric
properties by using a more sagacious constructing method of the
metric. As a partial negative solution to this problem, we prove the
following non-existence result:

\begin{prop}\label{propno}
Let $\Sigma^n$ be a homotopy sphere which admits a metric $g$ and an
isoparametric function $f$ with 2 points as the focal set. Suppose
that the level hypersurfaces are all totally umbilic. Then
$\Sigma^n$ is diffeomorphic to $S^n$.
\end{prop}

By combining Proposition \ref{propno} with a topological argument, we obtain

\begin{thm}\label{nonisopara}
Every odd dimensional exotic sphere admits no totally isoparametric
functions with 2 points as the focal set.
\end{thm}

Recall that a \emph{totally isoparametric} function is an
isoparametric function so that each regular level hypersurface has
constant principal curvatures, as defined in \cite{GTY14}. As is
well known, an isoparametric function on a unit sphere must be
totally isoparametric.

Last but not least, we have both existence and non-existence results
of isoparametric functions on some homopoty spheres which also have
$SC^p$-property ( see \cite{Be78} for definition ).

The paper is organized as follows. In Section 2, we review some
basic definitions and give a proof to Theorem \ref{Main}. Then we
apply the construction in Theorem \ref{Main} to homotopy spheres and
the Eells-Kuiper projective planes in Section 3. And in Section 4,
we show Proposition \ref{propno} and Theorem \ref{nonisopara}. The
last section will be concerned with the existence and non-existence
of isoparametric functions on exotic spheres with $SC^p$-property.
In particular, Theorem 5.1 improves a beautiful result of
B\'{e}rard-Bergery \cite{BB77}.

\section{A fundamental construction}
This section is devoted to a complete proof of Theorem \ref{Main}.
Let's start by recalling the definition of a Morse-Bott function. A
smooth function $f$ on a smooth manifold $N$ is called a
\emph{Morse-Bott function} or a generalized Morse function if: (i)
the critical set defined by equation $\mathrm{d}f=0$  consists of a
family of smooth submanifolds, the so called critical submanifolds;
(ii) the Hessian $H_f$ is a non-degenerate quadratic form in the
normal direction of each critical submanifold to $N$.

In differential topology, it is well known that the topology of a
manifold can be analyzed by studying Morse-Bott functions on it. The
following result is essentially due to Wang \cite{Wa87}.
\begin{prop}\label{transMorse}
A transnormal function $f$ on a complete connected Riemannian manifold
$N$ is a Morse-Bott function. Moreover, the critical set $C(f)$ of $f$
coincides exactly with the focal set of $f$.
\end{prop}
\begin{proof}
Combine Lemma 4 with Lemma 6 in \cite{Wa87}.
\end{proof}

Owning to Proposition \ref{transMorse}, it is natural to propose the following converse problem.
\begin{prob}
Let $f$ be a Morse-Bott function on a smooth manifold $N$. When does
$N$ carry a Riemannian metric $g$ so that $f$ is transnormal with
respect to this $g$?
\end{prob}
This problem motivates us to establish Theorem \ref{trans} below.
Let's first investigate the topology of a Morse-Bott function whose
critical set has only two components.

\begin{prop}
Let $N$ be a connected closed manifold and $f:N\rightarrow\mathbb{R}$ a Morse-Bott
function so that $C(f)$
 has two connected closed components, say $M_+$ and $M_-$. Then $N$ is diffeomorphic to
 the glued manifold
$$D(M_-)\sqcup M\sqcup_{\varphi} D(M_+),$$
for some diffeomorphism $\varphi:\partial D(M_-)\rightarrow\partial
D(M_+)$, where $M=\partial D(M_-)\times [0, 1]$, $D(M_+)$ and
$D(M_-)$ are respectively the normal disc bundles over $M_+$ and
$M_-$.
\end{prop}
\begin{proof}
By the generalized Morse lemma and regular interval theorem ( see
\cite{Hi76}, pp. 149, 153 ).
\end{proof}

\begin{rem}
Duan and Rees \cite{DR92} also considered the situation where the critical set of
the smooth map
$f:N\rightarrow\mathbb{R}$ has only two components, both of which
are smooth submanifolds of $N$. Actually, they did not assume that the critical set
 of $f$ are non-degenerate in any sense.
\end{rem}

Next, we set about constructing the desired metric. The construction is a modified version of A. Weinstein's
 construction ( see, for example, \cite{Be78}, pp. 231-233 ).
\begin{thm}\label{trans}
Let $N$ be a closed smooth manifold, and $f:N\rightarrow\mathbb{R}$
a Morse-Bott function with $C(f)=M_+\sqcup M_-$, where $M_+$ and
$M_-$ are both closed connected submanifolds of codimensions more
than $1$. Then there exists a Riemannian metric on $N$ so that $f$
 is transnormal. In fact, the metric can be chosen so that $M_+$ and
 $M_-$ are both totally geodesic.
\end{thm}
\begin{proof}
Without loss of generality, we can assume Im$f=[\alpha, \beta]$  and
$M_+=f^{-1}(\beta)$, $M_-=f^{-1}(\alpha)$. Our construction of the
metric is divided into three steps.

\setlength{\parindent}{0mm}
\textbf{Step 1:}
\setlength{\parindent}{2em}

Firstly, by the generalized Morse lemma, there exists a vector
bundle ( a
 tubular neighborhood )
 over $M_-$ with some metric, say $\xi^-=(\pi_-, E_-, M_-)$ and an embedding
 $\sigma_-:\mathcal{U}\rightarrow N$, where $\mathcal{U}\subset E_-$ is an open
 neighborhood of the zero section, such that

(1) $\sigma_-|_{M_-}=\mathrm{id}$, where $M_-$ is identified with the zero section
 of $\xi^-$ and $\sigma_-(\mathcal{U})$ is an open neighborhood of $M_-$ in $N$.
  Notice here rank $\xi^-=$dim $N-$dim $M_-\geq 2$.

(2) the composition $\mathcal{U}\xrightarrow{\sigma_-}N\xrightarrow{f}\mathbb{R}$
 is given by
$$f(\sigma_-(p, v))=|v|^2+\alpha$$
 for any $(p, v)\in \pi_-^{-1}(p)\cap\mathcal{U}$ and $p\in M_-$.

Now, since $M_-$ is compact, we can choose $\epsilon>0$ so small
that
$$M_-(\epsilon):=\{x\in N|~f(x)\leq\alpha+\epsilon^2\}\subset\sigma_-(\mathcal{U}).$$
Thus, we have a disc bundle,
$D_{\xi^-}(\epsilon):=\{(p,v)\in E_-|p\in M_-, |v|\leq\epsilon\}$
with $S_{\xi^-}(\epsilon):=\partial D_{\xi^-}(\epsilon)$ so that $\sigma_-(D_{\xi^-}(\epsilon))=M_-(\epsilon)$.

Next, for each dimension $d$, let us choose a smooth Riemannian
metric $g^d$ on the disc
$D^d(\epsilon):=\{x\in\mathbb{R}^d|~|x|\leq\epsilon\}$ as
 $g^d=F(r)^2\mathrm{d}r^2+G(r)^2\mathrm{d}s_{d-1}^2$, where
$\mathrm{d}s_{d-1}^2$ is the canonical metric on $S^{d-1}(1)$, and
$F, G$ are smooth functions on $[0, \epsilon]$ with
\[
F(r)=\left\{
\begin{array}{cc}
1, & \mathrm{for} ~r \leq \frac{\epsilon}{2}~;\\
2r, & \mathrm{for} ~r\geq \frac{3\epsilon}{4}~,
\end{array}\right.
\]
and
\[
G(r)=\left\{
\begin{array}{cc}
r,& \mathrm{for} ~r \leq \frac{\epsilon}{2}~;\\
1,& \mathrm{for} ~r \geq \frac{3\epsilon}{4}~.
\end{array}\right.
\]

The choice of the metric ensures the following two properties:

(1) $g^d$ is $O(d)$-invariant under the usual action of the
orthogonal group $O(d)$ with the origin $o\in\mathbb{R}^d$ fixed;

(2) $g^d$ is a product near the boundary $\partial D^d(\epsilon)=S^{d-1}(\epsilon)$.
Moreover, $|\frac{\partial}{\partial r}|=2r$ near the boundary, where $r=|x|$, and
 $\frac{\partial}{\partial r}$ is the radius direction.

With this metric, each radius has the same length, say $\delta$. Furthermore, it
is not difficult to show each radius is a geodesic with arc length as parameter in $D^d(\epsilon)$.

Let $h_-$ be any Riemannian metric on $M_-$. By a construction due
to J. Vilms \cite{Vi70}, there exists a Riemannian metric $h_-'$ on
$D_{\xi^-}(\epsilon)$, such that the fibration
$\pi:(D_{\xi^-}(\epsilon), h_-')\rightarrow(M_-, h_-)$ is a
Riemannian submersion with totally geodesic fibers isometric to
$(D^{d_-}(\epsilon), g^{d_-})$ with $d_-=$ rank $\xi^-$. In the case
at hand, this construction can be described as follows. Note that we
have chosen the Euclidean metric on $(\pi_-, E_-, M_-)$, and let us
choose a metric linear connection on $(\pi_-, E_-, M_-)$. This
connection produces a distribution $\mathcal{H}_-$ on the tangent
bundle $TE_-$ of $E_-$, complementary to the vertical distribution
$\mathcal{V}_-$ ( kernel of $\mathrm{d}\pi_-$ ), and the associated
parallel transport is by linear isometries. Then the desired metric
$h_-'$ on $D_{\xi^-}(\epsilon)$ can be constructed in the following
way: one takes the vertical and horizontal distributions to be
orthogonal, lifts $h_-$ from $M_-$ to the horizontal distribution
and then takes the metric $g^{d_-}$ on the fiber $D^{d_-}(\epsilon)$.
Clearly parallel transport is still an isometry for these metrics on
the fibers, so that the fibers are totally geodesic.

Notice that the metric $h_-'$ is still a product near the boundary
$S_{\xi^-}(\epsilon)$ of $D_{\xi^-}(\epsilon)$ and the zero section
is totally geodesic in $D_{\xi^-}(\epsilon)$. Moreover, any geodesic
issuing from a point in $S_{\xi^-}(\epsilon)$ and orthogonal to
$S_{\xi^-}(\epsilon)$ in $D_{\xi^-}(\epsilon)$ stays in the same
fiber of the fibration. In particular, all the geodesics orthogonal to
the boundary $S_{\xi^-}(\epsilon)$ in $D_{\xi^-}(\epsilon)$ have the
same length $2\delta$.

Now, by the diffeomorphism $\sigma_-:D_{\xi^-}(\epsilon)\rightarrow
M_-(\epsilon)$, we obtain an induced metric on $M_-(\epsilon)$, say
$h_-(\epsilon)$. On the other side, we can also construct an
analogous Riemannian metric, say $h_+(\epsilon)$, on
$M_+(\epsilon):=\{x\in N|~f(x)\geq \beta-\epsilon^2\}$ which is
diffeomorphic to a disc bundle $D_{\xi^+}(\epsilon)$ over $M_+$,
where $\epsilon$ is small enough, the same as before.

\setlength{\parindent}{0mm}
\textbf{Step 2:}
\setlength{\parindent}{2em}

The task we face is to extend smoothly the metrics $h_-(\epsilon)$
and $h_+(\epsilon)$ on $M_-(\epsilon)$
 and $M_+(\epsilon)$ to the whole manifold $N$ with appropriate properties.
We are left to work on $M:=\{x\in N|~\alpha+\epsilon^2\leq f(x) \leq
\beta-\epsilon^2\}$, since $N=M_-(\epsilon)\sqcup M\sqcup
M_+(\epsilon)$.

To do it, by regular interval theorem ( see \cite{Hi76} ), $\partial
M_-(\epsilon)$ is diffeomorphic to $\partial M_+(\epsilon)$. More
precisely, there is a diffeomorphism $\tau:M\rightarrow
S_{\xi^-}(\epsilon)\times[\alpha+\epsilon^2, \beta-\epsilon^2]$
satisfying
$$\tau(f^{-1}(a))=S_{\xi^-}(\epsilon)\times \{a\},~\forall a\in [\alpha+\epsilon^2, \beta-\epsilon^2].$$

Let us now choose a family $g_t$ of Riemannian metrics on
$S_{\xi^-}(\epsilon)$ , which are smooth in the parameter $t$ for
$t\in[\alpha+\epsilon^2, \beta-\epsilon^2]$ and constant near
$\alpha+\epsilon^2$ and $\beta-\epsilon^2$, respectively. Moreover,
we ask that $g_{\alpha+\epsilon^2}$ is the metric induced from
$h_-(\epsilon)$, and $g_{\beta-\epsilon^2}$ is induced from
$h_+(\epsilon)$ by the diffeomorphism $\tau$.

Then a metric $h_0$ on $S_{\xi^-}(\epsilon)\times [\alpha+\epsilon^2, \beta-\epsilon^2]$
 can be constructed in such a way that the direct decomposition given by the projections
 onto the two factors is orthogonal, the metric induced on each
 $S_{\xi^-}(\epsilon)\times \{t\}$ is $g_t$, and the metric induced on
 $\{p\}\times [\alpha+\epsilon^2, \beta-\epsilon^2]$ is the canonical metric
 with length  $\beta-\alpha-2\epsilon^2$ on $[\alpha+\epsilon^2, \beta-\epsilon^2]$
 for each point $p\in S_{\xi^-}(\epsilon)$. It is not difficult to show that for any given
  point $p\in S_{\xi^-}(\epsilon)$, the curve
$$\gamma: [\alpha+\epsilon^2, \beta-\epsilon^2]\rightarrow S_{\xi^-}(\epsilon)\times [\alpha+\epsilon^2, \beta-\epsilon^2],~t\mapsto(p,t),$$
is a geodesic in $S_{\xi^-}(\epsilon)\times [\alpha+\epsilon^2, \beta-\epsilon^2]$.

\setlength{\parindent}{0mm}
\textbf{Step 3:}
\setlength{\parindent}{2em}

Lastly, these metrics $h_-(\epsilon)$, $\tau^*h_0$ and
$h_+(\epsilon)$ on $M_-(\epsilon)$, $M$ and $M_+(\epsilon)$
respectively can be glued together to give a desired metric $g_N$ on
$N$. Observe that $|\nabla f|\equiv1$ on $M$. On the other hand,
near the boundary of $M_-(\epsilon)$, one has $|\nabla f|=|2r\nabla
r|=|\frac{1}{2r}\frac{\partial}{\partial r}|=1$, which also holds
near the boundary of $M_+(\epsilon)$. These properties guarantee
that the glued metric $g_N$ on $N$ is smooth. From this
construction, it is obvious that every geodesic starting
orthogonally from any point in $M_-$($M_+$) must arrive at
$M_+$($M_-$) orthogonally with length
$2\delta+\beta-\alpha-2\epsilon^2$.

By a direct computation, we see
\[
|\nabla f|^2=\left\{
\begin{array}{cc}
\frac{4(f-\alpha)}{F^2(\sqrt{f-\alpha})}& \alpha \leq f\leq \alpha+\epsilon^2 \\
1& \alpha+\epsilon^2\leq f \leq \beta-\epsilon^2\\
\frac{4(\beta-f)}{F^2(\sqrt{\beta-f})}& \beta-\epsilon^2 \leq f\leq \beta ,
\end{array}\right.
\]
which means that $f$ is transnormal. Now, the proof is complete.
\end{proof}
\begin{rem}
In general, $M_+$ and $M_-$ of a transnormal function $f$ are not necessarily totally geodesic.
\end{rem}

According to the proof of Theorem \ref{trans}, we have the following enjoyable observation.
\begin{prop}
The function $f$ is isoparametric on
the domain $N-M$ with respect to the metric constructed in Theorem \ref{trans}.
\end{prop}
\begin{proof}
We only need to verify $f$ on
$M_-(\epsilon)=\sigma_-(D_{\xi^-}(\epsilon))$. Since we are working
on a Riemannian submersion with totally geodesic fibers, it follows
that $\triangle$ on $M_-(\epsilon)$ has a decomposition
$\triangle=\triangle^{\mathcal{H}}+\triangle^{\mathcal{V}}$, where
$\triangle^{\mathcal{H}}$ and $\triangle^{\mathcal{V}}$ are the
horizontal and vertical Laplacians on $M_-(\epsilon)$ with respect
to the Riemannian submersion $\pi\circ
\sigma_-^{-1}:M_-(\epsilon)\rightarrow M_-$. In our case, it is not
difficult to see that $\triangle^{\mathcal{H}}f=0$. Thus, by a
direct computation,
$$\triangle f=\triangle^{\mathcal{V}}(r^2)=\frac{2}{F^2}-\frac{2rF_r}{F^3}+\frac{2(\mathrm{dim}N-\mathrm{dim}M_--1)rG_r}{F^2G},$$
that is, $\triangle f$ is a continuous function of $f$ on $M_-(\epsilon)$.
The proof then follows.
\end{proof}

Up to now, we have constructed a metric $g_N$ on $N$ such that $f$
is transnormal on the whole $N$, and isoparametric near $M_-$ and $M_+$.
In order to prove Theorem \ref{Main}, we have to deform the metric
$g_N$ on $N$ so that $f$ is isoparametric. For this purpose, we
prepare three lemmas.

\begin{lem}[Local version of Moser's theorem]\label{local}
Let $Q^n=(0,1)^n$ and $\bar{Q}^n=[0,1]^n$ be the open and closed
unit cubes in Euclidean space, respectively. Let $f$, $g$ be two
positive smooth functions on $\bar{Q}^n$ such that $f=g$ near
$\partial \bar{Q}^n$ and $\int_{Q^n}(f-g)\mathrm{d}x=0$. Then there
exists a diffeomorphism $\psi:Q^n\rightarrow Q^n$ satisfying

(1)~ $g\circ\psi~\mathrm{det}\nabla\psi=f  \;~in~Q^n;$

(2)~ $\psi=\mathrm{id}$ near $\partial \bar{Q}^n$;

(3)~ $\psi$ is isotopic to $\mathrm{id}$.

In fact, the diffeomorphism $\psi$ can be decomposed as
$\psi=\psi_n\circ\psi_{n-1}\circ...\circ\psi_1:Q^n\rightarrow Q^n$,
and for any $1\leq s\leq n$,
$$\psi_s:Q^n\rightarrow Q^n, (x_1,x_2,...,x_n)\mapsto(x_1, x_2,...,x_{s-1},\upsilon_s(x_1, x_2,...,x_n),x_{s+1},...,x_n)$$
is a diffeomorphism such that $\frac{\partial \upsilon_s}{\partial
x_s}>0$, $\psi_s=\mathrm{id}$ near $\partial \bar{Q}^n$, and
$\psi_s$ preserves the line segments $\{(x_1,x_2,...,x_n)\in Q^n|~
x_i=c_i, 1\leq i\leq n, i\neq s\}$.
\end{lem}

\begin{lem}[Global version of Moser's theorem]\label{global}
Let $\tau$ and $\sigma$ be two volume elements on a compact
connected manifold $M$ with
$$\int_{M}\tau=\int_{M}\sigma.$$
Then there exists a diffeomorphism $\psi: M\rightarrow M$ so that
$\psi$ is isotopic to $\mathrm{id}:M\rightarrow M$ and satisfies
$\psi^*\tau=\sigma$.
\end{lem}

\begin{rem}
For the proofs of Lemma \ref{local} and Lemma \ref{global}, we refer
to \cite{Mo65} and \cite{DM90}. However, we should indicate that
Lemma 2 was not proved correctly in \cite{Mo65}, and it was
rectified in \cite{DM90}. As Calabi pointed out that the global
version of Moser's theorem holds even for non-orientable manifolds
if one uses the concept of ``odd'' forms ( see \cite{Mo65} ).
\end{rem}

At last, we need a lemma for deforming metrics.
\begin{lem}\label{conformal}
Let $\mathcal{S}^{n-1}$ be an $(n-1)$-dimensional smooth manifold
and $M^n=\mathcal{S}^{n-1}\times [\alpha_0, \beta_0]$. Let $h$ be a
metric on $M$ defined by a smooth one-parameter family of metrics
$\{h_t| ~t\in [\alpha_0, \beta_0]\}$ on $\mathcal{S}$, i.e.
$h|_{(x,t)}=h_t+\mathrm{d}t^2$ for $(x,t)\in \mathcal{S}^{n-1}\times
[\alpha_0, \beta_0]$. Define a new metric $\tilde{h}$ on $M$ given
by $\{\tilde{h}_t|~ \tilde{h}_t=e^{2u}h_t,~t\in[\alpha_0,
\beta_0]\}$ depending only on a smooth function $u$ on $M$.
Considering the smooth projection function $f:
\mathcal{S}^{n-1}\times [\alpha_0, \beta_0]\rightarrow \mathbb{R},
(x,t)\mapsto t$, we have
$$\tilde{\triangle}f=\triangle f+(n-1)\frac{\partial u}{\partial t},$$
where $\tilde{\triangle}$ and $\triangle$ are the Laplacians with respect to $\tilde{h}$ and $h$ respectively.
\end{lem}
\begin{proof}
By a straightforward computation.
\end{proof}

After the long preparation, we are ready to prove theorem \ref{Main}.

\setlength{\parindent}{0mm}\textbf{Proof of Theorem \ref{Main}:}
\setlength{\parindent}{2em} Assume the same notations as in Theorem
\ref{trans}. By the proof of Theorem \ref{trans}, there exists a
diffeomorphism $\varphi: S_{\xi^-}(\epsilon)=\partial
D_{\xi^-}(\epsilon)\rightarrow S_{\xi^+}(\epsilon)=\partial
D_{\xi^+}(\epsilon)$ such that $N$ is diffeomorphic to the glued
manifold $D_{\xi^-}(\epsilon)\sqcup M'\sqcup_{\varphi}
D_{\xi^+}(\epsilon)$, where $M'=S_{\xi^-}(\epsilon)\times
[\alpha+\epsilon^2, \beta-\epsilon^2] $. Additionally, we have
constructed metrics $h_-'$, $h_+'$ and $h_0$ on
$D_{\xi^-}(\epsilon)$, $D_{\xi^+}(\epsilon)$ and $M'$ respectively.
As we asserted before, $f$ is isoparametric on $D_{\xi^-}(\epsilon)$
and $D_{\xi^+}(\epsilon)$.

Let $\omega_-$ and $\omega_+$ be the volume elements of
$(S_{\xi^-}(\epsilon), h_-'|_{S_{\xi^-}(\epsilon)})$ and
$(S_{\xi^+}(\epsilon), h_+'|_{S_{\xi^+}(\epsilon)})$ respectively.
Then $(S_{\xi^-}(\epsilon),
g_{\beta-\epsilon^2})=(S_{\xi^-}(\epsilon),
\varphi^*(h_+'|_{S_{\xi^+}(\epsilon)}))$ has the volume element
$\varphi^*\omega_+$. By Lemma \ref{global} ( Global version of
Moser's theorem ), there exists a diffeomorphism
$\psi:S_{\xi^-}(\epsilon)\rightarrow S_{\xi^-}(\epsilon)$, which is
isotopic to $\mathrm{id}$ and satisfies
$$\psi^*(\varphi^*\omega_+)=\lambda \omega_-,$$
where $\lambda=\frac{\int \varphi^*\omega+}{\int \omega_-}$.
Defining $\tilde{\varphi}:=\varphi\circ\psi$, we see that
$\tilde{\varphi}$ is isotopic to $\varphi$ and
$$\tilde{\omega}=(\varphi\circ\psi)^*\omega_+=\psi^*(\varphi^*\omega_+)=\lambda\omega_-,$$
where $\tilde{\omega}$ is the volume element of $(S_{\xi^-}(\epsilon), \psi^*g_{\beta-\epsilon^2})$.
Since $\tilde{\varphi}$ is isotopic to $\varphi$, it follows that
$D_{\xi^-}(\epsilon)\sqcup M'\sqcup_{\tilde{\varphi}} D_{\xi^+}(\epsilon)$ is diffeomorphic
to $D_{\xi^-}(\epsilon)\sqcup M'\sqcup_{\varphi} D_{\xi^+}(\epsilon).$

Consider now $\tilde{N}:=D_{\xi^-}(\epsilon)\sqcup
M'\sqcup_{\tilde{\varphi}} D_{\xi^+}(\epsilon)$. Let
$\tilde{f}:\tilde{N}\rightarrow\mathbb{R}$ be a function such that
$\tilde{f}|_{D_{\xi^{\pm}}(\epsilon)}=f\circ\sigma_{\pm}|_{D_{\xi^{\pm}
}(\epsilon)}$ and $\tilde{f}|_{M'}=f\circ\tau^{-1}|_{M'}$. It is
obvious that $\tilde{N}$ is diffeomorphic to $N$ and $\tilde{f}$ is
a Morse-Bott function with the critical submanifolds
$C(\tilde{f})=\tilde{M}_+\bigsqcup \tilde{M}_-$, where $\tilde{M}_+$
and $\tilde{M}_-$ are diffeomorphic to $M_+$ and $M_-$,
respectively. With respect to the metrics constructed in Theorem
\ref{trans}, $\tilde{f}$ is isoparametric on $(D_{\xi^-}(\epsilon),
h_-')$ and $(D_{\xi^+}(\epsilon), h_+')$. On $M'$, choose a metric
$g$ defined by a smooth family of metrics $\{g_t'|~t\in
[\alpha+\epsilon^2, \beta-\epsilon^2]\}$ on $S_{\xi^-}(\epsilon)$
such that
$$g_{\alpha+\epsilon^2}'=g_{\alpha+\epsilon^2},~g_{\beta-\epsilon^2}'=\psi^*g_{\beta-\epsilon^2}$$
 and $g_t'$ is
constant near $\alpha+\epsilon^2$ and $\beta-\epsilon^2$.
 By lemma \ref{conformal},
let $u:M'\rightarrow\mathbb{R}$ be a smooth function to be determined and $\tilde{g}$
a new metric on $M'$ determined by the new family of metrics
$\{\tilde{g}_t'|~\tilde{g}_t'=e^{2u}g_t', t\in [\alpha+\epsilon^2, \beta-\epsilon^2]\}$,
then
$$\tilde{\triangle}\tilde{f}=\triangle \tilde{f}+(n-1)\frac{\partial u}{\partial t},$$
where $n=$ dim $N$, $\tilde{\triangle}$ and $\triangle$ are
Laplacians with respect to $\tilde{g}$ and $g$, respectively.

Next let $h:[\alpha+\epsilon^2,
\beta-\epsilon^2]\rightarrow\mathbb{R}$ be a smooth function with
 $h(t)=\triangle \tilde{f}|_{S_{\xi^-}(\epsilon)\times\{t\}}$ for $t$ near $\alpha+\epsilon^2$ and $\beta-\epsilon^2$.
 In order to find such an $u$ that $\tilde{f}$ is isoparametric on $(M', \tilde{g})$, it is sufficient to construct $u$
 as the solution
 of the following equation
$$\tilde{\triangle}\tilde{f}=h$$
or equivalently
$$h(t)=\triangle \tilde{f}+(n-1)\frac{\partial u}{\partial t}$$
which implies that for any $(x,t) \in M'$
$$(n-1)(u(x,t)-u(x,\alpha+\epsilon^2))=\int_{\alpha+\epsilon^2}^{t}[h(s)-\triangle \tilde{f}(x,s)]\mathrm{d}s.$$

On the other hand,
 by a direct computation,
 we have $\triangle \tilde{f}=\frac{1}{2}\frac{\partial}{\partial t}(\mathrm{log}~\mathrm{det}(g_{ij}))$,
where $g_t'=\sum_{i,j=1}^{n}g_{ij}(x,t)\mathrm{d}x_i\mathrm{d}x_j$
with respect to local coordinates $(x_1, x_2,..., x_{n-1})$ around
$x\in S_{\xi^-}(\epsilon)$. It follows that
$$\int_{\alpha+\epsilon^2}^{\beta-\epsilon^2}\triangle\tilde{f}(x,s)\mathrm{d}s
=\mathrm{log}\frac{\sqrt{\mathrm{det}(g_{ij}(x,~\beta-\epsilon^2))}}
{\sqrt{\mathrm{det}(g_{ij}(x,~\alpha+\epsilon^2))}}$$
which is just
equal to $\mathrm{log}\lambda$, independent of the choice of $x\in
S_{\xi^-}(\epsilon)$. Hence, an appropriate choice of the function
$h$ satisfies the equation
$\int_{\alpha+\epsilon^2}^{\beta-\epsilon^2}[h(s)-\triangle
\tilde{f}(\cdot, s)]\mathrm{d}s\equiv0$. At last
 we get a solution
$$u(x,t)=\frac{1}{n-1}\int_{\alpha+\epsilon^2}^{t}[h(s)-\triangle \tilde{f}(x,s)]\mathrm{d}s,$$
for $(x,t)\in M'$, such that $u|_{S_{\xi^-}(\epsilon)\times
\{t\}}\equiv0$ for $t$ near $\alpha+\epsilon^2$ and
$\beta-\epsilon^2$. The argument above implies that $\tilde{g}$
still has the product metric near the boundary of $M'$ as before.
 Therefore the three parts $(D_{\xi_-}(\epsilon), h_-')$, $(D_{\xi_+}(\epsilon), h_+')$ and $(M', \tilde{g})$ can be glued to gain a smooth metric $\tilde{g}_{\tilde{N}}$ on $\tilde{N}$
 so that $\tilde{f}$ is isoparametric.

The proof is now complete.$\hfill{\square}$

\section{Applications and examples}
In this section, the general construction in Theorem \ref{Main} is
applied to the existence of isoparametric functions on various
exotic spheres and Eells-Kuiper projective planes.

First, we give a short

\setlength{\parindent}{0mm}
\textbf{Proof of Corollary \ref{homosphere}:}
\setlength{\parindent}{2em}

\begin{proof}
A remarkable result of S. Smale ( see \cite{Sm61} and \cite{Mil07},
p. 128 ) states that there exists a Morse function with 2 critical
points
  on each homotopy sphere with dimension more than 4.
  Thus, it follows from Theorem \ref{Main} that each homotopy sphere with dimension more than 4 admits a metric and an isoparametric function with 2 points as the focal set.

The proof is now complete.
\end{proof}
\begin{rem}
As Smale's result is a breakthrough in differential topology,
Corollary \ref{homosphere} may be regarded as its continuation in
differential geometry.
\end{rem}

More generally, we also investigate the existence
 of isoparametric functions on exotic spheres such that at least one component of the critical set is not a single point.

Following J. Milnor,
 let $g:S^{m}\times S^{n}\rightarrow S^{m}\times S^{n}$ be
  an orientation-preserving diffeomorphism.
   Then a smooth manifold $M^{m+n+1}(g)$ depending on $g$ is obtained from
   disjoint spaces $\mathbb{R}^{m+1}\times S^{n}$ and $S^m\times\mathbb{R}^{n+1}$ by matching $(tx, y)$ and $(x', t'y')$, where $(x', y')=g(x, y)$ and $t'=\frac{1}{t}$ for any $x\in S^{m}$, $y\in S^n$ and $t\in (0, +\infty)$.
     In light of this construction, we have the following observation.
\begin{prop}\label{Morse-Bott}
Define a function $f$ on $M^{m+n+1}(g)$ by $$f:\mathbb{R}^{m+1}\times S^{n}\rightarrow\mathbb{R}, (X, y)\mapsto\frac{|X|^2}{1+|X|^2},$$
and
$$f:S^m \times \mathbb{R}^{n+1}\rightarrow\mathbb{R}, (x, Y)\mapsto\frac{1}{1+|Y|^2}.$$
$f$ is a Morse-Bott function with $C(f)=\{0\}\times S^n\sqcup S^m\times\{0\}$. Moreover, $M^{m+n+1}(g)$ admits an isoparametric function with focal set diffeomorphic to $S^m\sqcup S^n$ under some metric.
\end{prop}
\begin{proof}
It is obvious that $f$ is a well defined smooth function.
Moreover, $f$ is a Morse-Bott function with critical set $\{0\}\times S^n\sqcup S^m\times\{0\}$.
Therefore, $M^{m+n+1}(g)$ admits a metric and an isoparametric function with focal set diffeomorphic to $S^m\sqcup S^n$ under some metric,
by virtue of Theorem \ref{Main}.
\end{proof}

In order to construct concrete examples, we have the following method by J. Milnor \cite{Mil59}.
 Starting with two smooth maps $f_1:S^m\rightarrow SO(n+1)$ and $f_2:S^n\rightarrow SO(m+1)$, a diffeomorphism $g: S^m\times S^n\rightarrow S^m\times S^n$ is defined by
 $$g(x,y):=(f_2^{-1}(f_1(x)\cdot y)\cdot x,f_1(x)\cdot y)$$
  for $x\in S^m$ and $y\in S^n$. Denote $M(f_1, f_2):=M(g)$. There is also another way to construct manifolds, the so called Milnor pairing. The Milnor pairing is a bilinear pairing given by
$$\beta_{m,n}:\pi_m SO(n)\otimes \pi_n SO(m)\rightarrow \pi_0 \mathrm{Diff^+}~S^{m+n}\rightarrow \Gamma_{m+n+1}=\Theta_{m+n+1},$$
where $\mathrm{Diff^+}S^{m+n}$ is the group of all orientation
preserving-diffeomorphisms from $S^{m+n}$ onto itself; $\Gamma_l$
denotes the group of oriented diffeomorphism classes of twisted
$l$-spheres, $\Theta_l$ denotes the group of oriented homotopy
$l$-spheres up to relation $h$-cobordent, and in fact they are
isomorphic for each $l$ ( \cite{Mil07}, p. 5 ). We remark that two
homotopy $l$-spheres, $l>4$, are $h$-cobordant if and only if they
are diffeomorphic ( ref. \cite{KM63} ).
 The following theorem of J. Milnor tells us how the two
 constructing
methods above are related.
\begin{thm}[\cite{Mil07}, p. 218]\label{Milnor}
Given two smooth maps $f_1: S^m \rightarrow SO(n)\hookrightarrow SO(n+1)$ and $f_2: S^n \rightarrow SO(m)\hookrightarrow SO(m+1)$, the
resulting manifold
$M(f_1, f_2)$ is diffeomorphic to $\beta_{m,~n}(f_1, f_2)$.
\end{thm}

Applying Proposition \ref{Morse-Bott} and Theorem \ref{Milnor} will
lead to the following examples.
\begin{examp}
Each homotopy 7-sphere admits a metric and an isoparametric function
with focal set dffeomorphic to $S^3\sqcup S^3$, since $\Theta_7=$ Im
$\beta_{3,~3}$ ( see \cite{Mil07}, p. 187 ).
\end{examp}
\begin{rem}
This example generalizes Proposition 4.1 in \cite{GT13}, where only
Milnor $7$-spheres ( those can be represented as $S^3$ bundles over
$S^4$ ) are considered. In fact, it is well known that there exists 14 exotic
spheres in dimension 7, and four of which are not Milnor spheres.
\end{rem}

Recall that the group of homotopy spheres $\Theta_n$ has an
important subgroup $bP_{n+1}$ defined by Kervaire and Milnor in
\cite{KM63}. A homotopy $n$-sphere $\Sigma^n$ represents an element
of $bP^{n+1}$ if and only if $\Sigma^n$ bounds a parallelizable
manifold. To determine the group $\Theta_n$, they did a systematic
research of such subgroups. In particular, they showed that
$\Theta_7=bP_8$ and $bP_{4m+2}=\mathbb{Z}_2$ or $0$. As in the
literature, a Kervaire sphere $\Sigma^{4m+1}$ is the boundary of the
plumbing of two copies of the unit disc bundle $D(TS^{2m+1})$ of the
tangent bundle of $S^{2m+1}$ and consequently it is a generator of
$bP_{4m+2}$ ( see \cite{KM63}, \cite{HH67} and \cite{MM80} ). For
such homotopy spheres, we have the following result.
\begin{examp}
For $m>0$, the Kervaire sphere $\Sigma^{4m+1}$ admits a metric and
an isoparametric function with focal set diffeomorphic to
$S^{2m}\sqcup S^{2m}$. The reason is that one can identify
$\Sigma^{4m+1}$ with $M(f, f)$, where $f: S^{2m}\rightarrow
SO(2m)\hookrightarrow SO(2m+1)$ is the characteristic map for the
tangent bundle of $S^{2m+1}$ ( see \cite{HH67} and \cite{MM80} ).
\end{examp}
Additionally, more examples concerning an even-dimensional homotopy
sphere or an odd-dimensional homotopy sphere which is not
the boundary of any parallelizable manifold are given.
\begin{examp}
According to D. L. Frank \cite{Fr68}, it is known that $\Theta_8=$
Im $\beta_{3,~4}$, $2\Theta_{10}=$ Im $\beta_{3,~6}$. Hence, every
homotopy $8$-sphere admits a metric and an isoparametric function
with focal set diffeomorpic to $S^3\sqcup S^4$; every element in
$2\Theta_{10}$ admits a metric and an isoparametric function with
focal set diffeomorphic to $S^3 \sqcup S^6$.

It is also proved by D. L. Frank in \cite{Fr68} that there exists a
homotopy $15$-sphere $\Sigma^{15}$ which is contained in Im
$\beta_{11,~3}$ but not in $bP_{16}$. Therefore, this $\Sigma^{15}$
admits an isoparametric function with focal set $S^{11}\sqcup S^3$
for a certain metric, and it is not the boundary of any parallelizable
manifold.
\end{examp}
In fact, Proposition \ref{Morse-Bott} gives us more examples than we
show here ( see \cite{Mil59} and related references ).
\begin{rem}
It is interesting to compare our existence results of isoparametric functions on homotopy
spheres with the complete classification of cohomogeneity one actions on them
( see \cite{St96} ). According to \cite{St96}, aside from cohomogeneity one actions
on the standard spheres, only exotic Kervaire spheres admit cohomogeneity one
actions. More precisely, every exotic Kervaire sphere can be represented as
a smooth submanifold $\Sigma^{2n-1}_d$ of $\mathbb{C}^{n+1}$
defined by the equations
$$z_0^d+z_1^2+\cdots+z_n^2=0,~|z_0|^2+|z_1|^2+\cdots+|z_n|^2=1,$$
for some odd number $n$ and $d$, the so called Brieskorn variety ( see \cite{Hi66} ).
What is better, as discovered in \cite{HH67},
the Brieskorn variety $\Sigma^{2n-1}_d$ has a cohomogeneity one action by $SO(2)SO(n)$, which
can be defined as
$$(e^{i\theta}, A)(z_0, z_1,\ldots, z_n):=(e^{2i\theta}z_0, e^{di\theta}A(z_1,\ldots, z_n)^t)$$
for $e^{i\theta}$ in $SO(2)$ and $A$ in $SO(n)$. In fact, these
induced actions are all the cohomogeneity one actions on exotic
Kervaire spheres. Hence, all the exotic spheres but Kervaire spheres
do not admit cohomogeneity one actions. However, our existence
result implies that each exotic sphere with dimension more than 4
admits some Riemannian metric and an isoparametric function.
\end{rem}

At last, we will investigate the existence of isoparametric
functions on Eells-Kuiper projective planes. Recall that a closed
manifold is called an \emph{Eells-Kuiper projective plane} if it
admits a Morse function with three critical points. Eells-Kuiper
[EK1] obtained many remarkable results on such manifolds. For
instance, they showed that the integral cohomology ring of an
Eells-Kuiper projective plane $M$ is isomorphic to that of the real,
complex, quaternionic or Cayley projective plane. Consequently, the
dimension of $M$ must be equal to $2,4,8$, or $16$. Moreover,
 for $m\in\{1, 4, 8\}$, an Eells-Kuiper projective plane $M^{2m}$ is
 diffeomorphic to $D(\xi)\sqcup_{\varphi} D^{2m}$ for some
 diffeomorphism $\varphi:S(\xi)=\partial D(\xi)\rightarrow S^{2m-1}$,
 where $D(\xi)$ and $S(\xi)$ are the associated disc bundle and sphere bundle for certain vector bundle of rank $m$ over $S^{m}$.
 It follows that on every Eells-Kuiper projective plane $M^{2m}$($m\neq 2$) we can construct a
 Morse-Bott function whose critical set has two components, a point and $S^m$. Therefore, by Theorem \ref{Main}, we get the following result.
\begin{prop}
For $m=4$ or $8$, each Eells-Kuiper projective plane $M^{2m}$ admits
a metric and an isoparametric function so that, one component of the
focal set is a single point and the other is diffeomorphic to $S^m$.
\end{prop}
\begin{rem}
For $m=2$, due to the possible existence of exotic homotopy 4-sphere,
 we cannot obtain the similar existence results for this case.
\end{rem}
\begin{rem}
On the projective planes $CP^2$, $HP^2$ and $CaP^2$ with
Fubini-Study metrics, there are isoparametric functions so that the
focal sets have two components, a single point, and $S^2$, $S^4$ and
$S^8$, respectively. In fact, these isoparametric foliations are all
homogeneous. For $CP^2$ and $HP^2$, one uses the Hopf fibrations to
derive the desired conclusion. As for $CaP^2$, one uses the
cohomogenity one action of $\mathrm{Spin(9)}$ on
$F_4/\mathrm{Spin(9)}$ which is isometric to $CaP^2$.
\end{rem}

\section{Non-existence results on exotic spheres}
In this section, Proposition \ref{propno} and Theorem
\ref{nonisopara} are proved. Although any exotic $n$-sphere with
$n>4$ admits an isoparametric function with 2 points as the focal
set by Theorem \ref{Main}, Proposition \ref{propno} means that
regular level hypersurfaces of such isoparametric functions on an
exotic sphere cannot have strong extrinsic geometric properties.
This is a way to distinguish exotic spheres from standard spheres by
geometry.

\setlength{\parindent}{0mm}\textbf{Proof of Proposition \ref{propno}:}
\setlength{\parindent}{2em}
\begin{proof}
Let $2\delta:=d(m_+, m_-)$, where $d$ is the distance function with
respect to the Riemannian metric $g$. The existence of thus
isoparametric function $f$ implies $\Sigma^n$ can be decomposed as
$$\Sigma^n=\mathrm{exp} D_{\delta}(T_{m_-}\Sigma^n)\cup\mathrm{exp} D_{\delta}(T_{m_+}\Sigma^n),$$
where
 $D_{\delta}(T_{m_-}\Sigma^n)$ and $D_{\delta}(T_{m_+}\Sigma^n)$ are the discs of radius $\delta$ at $T_{m_-}\Sigma^n$ and $T_{m_+}\Sigma^n$
 respectively. It follows that $\Sigma^n$ is diffeomorphic to $D^n(1)\sqcup_{\eta} D^n(1)$, where $\eta: S^{n-1}(1)\rightarrow S^{n-1}(1)$ is determined by the following diagram
\[
\begin{array}{ccc}
S^{n-1}(1)&\hookrightarrow  D^n(1) \xrightarrow{\rho_+^{\delta}}  D_{\delta}(T_{m_+}\Sigma^n)\xrightarrow{\mathrm{exp}_{m_+}}&\Sigma^n \\
\downarrow{\eta}&                                                               &\parallel\\
S^{n-1}(1)&\hookrightarrow  D^n(1) \xrightarrow{\rho_-^{\delta}}  D_{\delta}(T_{m_-}\Sigma^n)\xrightarrow{\mathrm{exp}_{m_-}}&\Sigma^n,
\end{array}
\]
where $\rho_+^{\delta}$ and $\rho_-^{\delta}$ are the canonical linear diffeomorphisms.

Since the regular hypersurfaces of $f$ are all totally umbilical, we can assume the regular hypersurface $X_t:=\{x\in \Sigma^n|~d(p, m_-)=t\}$ has constant principal curvature $\lambda(t)$ for any $t\in (0, 2\delta)$.
According to
 Kowalski and Vanhecke ( \cite{KV86}, Theorem 12 ), we have
$$\mathrm{exp}_{m_-}^*g=H(r)\sum_{i=1}^n(\mathrm{d}x_i)^2+\frac{1-H(r)}{r^2}(\sum_{i=1}^{n}x_i\mathrm{d}x_i)^2,$$
where
$H(r)=\mathrm{exp}\{2\int_0^{r}(\lambda(t)-\frac{1}{t})\mathrm{d}t\}$,
$r^2=\sum_{i=1}^{n}(x_i)^2\in (0, \delta^2]$, and $(x_1, x_2,...,
x_n)$ is the normal coordinate around $m_-$. Similarly, near $m_+$,
$$\mathrm{exp}_{m_+}^*g=H(r)\sum_{i=1}^n(\mathrm{d}y_i)^2+\frac{1-H(r)}{r^2}(\sum_{i=1}^{n}y_i\mathrm{d}y_i)^2,$$
where
$H(r)=\mathrm{exp}\{2\int_0^{r}(\lambda(2\delta-t)-\frac{1}{t})\mathrm{d}t\}$,
$r^2=\sum_{i=1}^{n}(y_i)^2\in (0, \delta^2]$, and $(y_1, y_2,...,
y_n)$ is the normal coordinate around $m_+$. By the diagram
introducing the diffeomorphism $\eta$
 and the formulae of the metric around $m_+$ and $m_-$, we conclude that $\eta: S^{n-1}(1)\rightarrow S^{n-1}(1)$ is an
 isometry.
  It follows at once that $\Sigma^n$ is diffeomorphic to $S^n(1)$, which completes the proof.
\end{proof}
At last, for odd dimensional homotopy spheres, we have Theorem
\ref{nonisopara} by a topological argument.

\setlength{\parindent}{0mm}\textbf{Proof of Theorem \ref{nonisopara}:}
\setlength{\parindent}{2em}

\begin{proof}
Let $\Sigma^{2n+1}$ be a homotopy sphere and $f:\Sigma^{2n+1}\rightarrow\mathbb{R}$ a
 totally isoparametric
  function with respect to a metric $g$ on $\Sigma^{2n+1}$. Then for any regular value $t$, the regular hypersurface $M_t:=f^{-1}(t)$
  has constant principal curvatures in $(\Sigma^{2n+1}, g)$. Let $m(t)$ be the number of distinct principal curvatures of $M_t$ and assume $m(t)>1$ for some $t$.
  Then the tangent bundle
   of $M_t$ has the decomposition as follows,
$$TM_t=T_1\oplus T_2\oplus\cdot\cdot\cdot\oplus T_{m(t)},$$
where  $T_1$, $T_2$,..., $T_{m(t)}$ are the principal distributions.
On the other hand, we observe that each regular hypersurface $M_t$
is diffeomorphic to $S^{2n}$, since $M_t$ is a geodesic hypersphere
of each of two focal points. Hence, the Euler characteristic classes
satisfy
$$0\neq e(TM_t)=e(T_1)e(T_2)\cdot\cdot\cdot e(T_{m(t)}).$$
However, the right side of the equality is equal to $0$, since
$\mathrm{rank} T_k<2n$,
each class $e(T_k)\in$ $H^{\mathrm{rank} T_k}(M_t)=H^{\mathrm{rank} T_k}(S^{2n})$ vanishes, for $k \in \{1, 2,..., m(t)\}$.
That is a contradiction,
which implies that $m(t)=1$ and the regular hypersurface $M_t$ is totally umbilical for any regular value $t$. The conclusion follows immediately from Proposition \ref{propno}.
\end{proof}
\begin{rem}
According to \cite{HH67} and \cite{St96}, there exists at leat one exotic
Kervaire sphere $\Sigma^{4m+1}$ which has a cohomogeneity one action.
Consequently,
by Proposition 2.3 in \cite{GT13}, $\Sigma^{4m+1}$ admits a
totally isoparametric function $f$ under an invariant metric.
However, each component of the focal set of $f$ is not just a point,
but a smooth submanifold. Hence, the assumption on the focal set in
Theorem \ref{nonisopara} is essential.
\end{rem}
Due to Theorem \ref{nonisopara}, it is reasonable to ask
\begin{prob}
Does there exist an even dimensional exotic sphere $\Sigma^{2n}$ $(n>2)$ which admits a metric and a totally isoparametric function with 2 points as the focal set?
\end{prob}

\section{Isoparametric functions and $SC^p$-property}
In this section, both isoparametric functions and $SC^p$-structure
are investigated simultaneously. First, we recall the  definition.
Let $(M,g)$ be a Riemannian manifold and $p$ be a point in $M$. If
all the geodesics issued from $p$ are simply closed geodesics with
the same length, $(M,g)$ is said to have the \emph{$SC^p$-property}
at $p\in M$.

For background knowledge and a systematic research, we refer
to the classic book \cite{Be78}.

Next, we give the following existence theorem which improves a beautiful
result of B\'{e}rard-Bergery \cite{BB77} ( see also p. 159,
pp. 234-235 in \cite{Be78} ).
\begin{thm}\label{SC}
On any homotopy sphere $\Sigma^{2n}$ in $2\Theta_{2n}(n\geq3)$,
there exists a Riemannian metric so that it possesses
$SC^p$-property at two points, say $m_+$ and $m_-$. Furthermore
under the same metric, there exists an isoparametric function $f$ on
$\Sigma^{2n}$ with focal set $C(f)=\{m_+, m_-\}$. The later property
means that
 $\Sigma^{2n}$ is locally harmonic at both points ${m_+}$ and ${m_-}$.
\end{thm}
\begin{proof}
We first observe that,
given any homotopy sphere $\Sigma^{2n}\in 2\Theta_{2n}$,
there exists an orientation-preserving diffeomorphism $\eta:S^{2n-1}\rightarrow S^{2n-1}$ such that
$$\Sigma^{2n}=D^{2n}\sqcup_{\eta\circ\eta}D^{2n}.$$
Let $(S^{2n-1}, g_0)$ be the standard unit sphere with volume element $\omega_0$ and the
standard antipodal map $\tau_0$.
By the constructions in Theorem \ref{trans} and Theorem \ref{Main},  to prove Theorem \ref{SC}, it suffices to
find an orientation-preserving diffeomorphism $\varphi: S^{2n-1}\rightarrow S^{2n-1}$ such that

(1) $\varphi$ is isotopic to $\eta\circ\eta$;

(2) $\varphi^*\omega_0=\omega_0$;

(3) $\varphi\circ\tau_0=\tau_0\circ\varphi$.

\setlength{\parindent}{0mm}
 The conditions (1) and (2) imply the existence of a metric $g$ and an isoparametric function $f$ with $C(f)=\{m_+, m_-\}$,
 and condition (3) guarantees that the metric $g$ satisfies the $SC^p$-property at $m_+$ and $m_-$.
\setlength{\parindent}{2em}

To construct the desired $\varphi$, let
$S_+:=\{x=(x_1,x_2,...,x_{2n})\in S^{2n-1}|~x_1\geq0\}$,
 and $S_-:=\{x=(x_1,x_2,...,x_{2n})\in S^{2n-1}|~x_1\leq0\}$. First, choose a diffeomorpism $\zeta:S^{2n-1}\rightarrow S^{2n-1}$,
 which is the restriction of an orientation-preserving diffeomorphism
  from $D^{2n}$ to $D^{2n}$,
  such that the composition $\eta':=\eta\circ\zeta$ satisfies $\eta'|_{S_+}=\mathrm{id}$ and $\eta'$ is isotopic to $\eta$.
  Thus we get two forms $(\eta')^*\omega_0$ and $\omega_0$ on $S^{2n-1}$, which are equal on the
  domain $S_+$ and have equal integral over the domain $S_-$. Applying Lemma \ref{local} to
  the two
  forms, we have a diffeomorphism $u:S^{2n-1}\rightarrow S^{2n-1}$ such that $u$ is
   isotopic
   to $\mathrm{id}$, $u|
_{S_+}=\mathrm{id}$, and $u^*((\eta')^*\omega_0)=\omega_0$. Define a
diffeomorphism $\tilde{\eta}:=\eta'\circ u$. It is clear that
$\tilde{\eta}|_{S_+}=\mathrm{id}$, $(\tilde{\eta})^*\omega_0=\omega_0$.

Now, the diffeomorpism $\varphi$ is defined by the compositions
$\varphi:=\tau_0\circ\tilde{\eta}\circ\tau_0\circ\tilde{\eta}$. We
are left to verify that $\varphi$ possesses the required properties.
For condition (1), we note that $\tau_0$ is isotopic to
$\mathrm{id}$ on the odd dimensional sphere $S^{2n-1}$. While
condition (2) is evident to verify. As for condition (3), since
$\tau_0\circ\tilde{\eta}\circ\tau_0$ is the identity on $S_-$, it
follows that $\tau_0\circ\tilde{\eta}\circ\tau_0$ commutes with
$\tilde{\eta}$, and thus $\varphi\circ\tau_0=\tau_0\circ\varphi$. It
completes the proof.
\end{proof}

\begin{rem}
The first half conclusion of Theorem 5.1 was proved by
B\'{e}rard-Bergery \cite{BB77} ( p. 237. see also p. 4, pp. 234-235
of \cite{Be78}). The later half means that every geodesic
hypersphere centered at $m_+$ or $m_-$ is of constant mean
curvature. In particular, there exists at least one exotic sphere of
dimension $10$ which possesses these properties.

\end{rem}

As for odd dimensional homotopy spheres, we have the following non-existence result.
\begin{thm}\label{noSC}
Let $\Sigma^7$ be a 7-dimensional homotopy sphere and
$f:\Sigma^7\rightarrow\mathbb{R}$ a transnormal function  with the
focal set $C(f)=\{m_+, m_-\}$ under some metric $g$. Suppose that
$(\Sigma^7, g)$ satisfies the $SC^p$-property at either $m_+$ or
$m_-$. Then $\Sigma^7$ is diffeomorphic to $S^7$ or the element
$[14]$ in $ \Theta_7 =\mathbb{Z}_{28}$.
\end{thm}
\begin{proof}
As in the proof of Proposition \ref{propno},
let $2\delta:=d(m_+, m_-)$, where $d$ is the distance function
with respect to the Riemannian metric $g$.
Since $f$ is a transnormal function with $C(f)=\{m_+, m_-\}$, $\Sigma^7$ has the following decomposition,
$$\Sigma^n=\mathrm{exp} D_{\delta}(T_{m_-}\Sigma^n)\cup\mathrm{exp} D_{\delta}(T_{m_+}\Sigma^n),$$
where $D_{\delta}(T_{m_-}\Sigma^7)$ and $D_{\delta}(T_{m_+}\Sigma^7)$ are the discs of radius $\delta$ at $T_{m_-}\Sigma^7$ and $T_{m_+}\Sigma^7$ respectively. What is more, $\Sigma^7$ is diffeomorphic
to $D^7(1)\sqcup_{\sigma} D^7(1)$, where $\sigma: S^{6}(1)\rightarrow S^{6}(1)$ is determined by the following diagram,
\[
\begin{array}{ccc}
S^{6}(1)&\hookrightarrow  D^7(1) \xrightarrow{\rho_+^{\delta}}  D_{\delta}(T_{m_+}\Sigma^7)\xrightarrow{\mathrm{exp}_{m_+}}&\Sigma^7 \\
\downarrow{\sigma}&                                                               &\parallel\\
S^{6}(1)&\hookrightarrow  D^7(1) \xrightarrow{\rho_-^{\delta}}  D_{\delta}(T_{m_-}\Sigma^7)\xrightarrow{\mathrm{exp}_{m_-}}&\Sigma^7,
\end{array}
\]
where $\rho_+^{\delta}$ and $\rho_-^{\delta}$ are the canonical linear diffeomorphisms.

The assumption that $(\Sigma^7, g)$ satisfies the $SC^p$-property at the point $m_+$($m_-$) implies that it satisfies the $SC^p$-property at the point $m_-$($m_+$).
That is to say, the assumption can deduce that $(\Sigma^7, g)$ satisfies the $SC^p$-property at both of the 2 points $m_+$ and $m_-$.
 Thus we obtain
$$\sigma\circ\tau_0=\tau_0\circ\sigma.$$

Define a map
\[
\begin{array}{ccc}
\mathcal{A}&:&\Theta_7=\mathbb{Z}_{28}\rightarrow\Theta_7=\mathbb{Z}_{28}.\\
 &&D^7\sqcup_{\varphi}D^7\mapsto D^7\sqcup_{\tau_0\circ\varphi\circ\tau_0}D^7
\end{array}
\]
It is obvious that $\mathcal{A}$ is a well-defined homomorphism. According to \cite{ADPR}, $\mathcal{A}$ is determined as follows,
\[
\begin{array}{ccc}
\mathcal{A}&:&\Theta_7=\mathbb{Z}_{28}\rightarrow\Theta_7=\mathbb{Z}_{28}.\\
&&[n]\mapsto[-n]
\end{array}
\]
Now, we have the diffeomorphsims
 $\mathcal{A}(\Sigma^7)=\mathcal{A}(D^7\sqcup_{\sigma}D^7)=D^7\sqcup_{\tau_0\circ\sigma\circ\tau_0}D^7=D^7\sqcup_{\sigma}D^7=\Sigma^7$.
 Hence, it follows that $\Sigma^7=[0]$ or $[14]\in \mathbb{Z}_{28}=\Theta_7$.
\end{proof}
\begin{rem}
Similarly, let $\Sigma^{15}$ be a 15-dimensional homotopy sphere which bounds a parallelizable manifold and $f:\Sigma^{15}\rightarrow\mathbb{R}$ a transnormal function  with the focal set $C(f)=\{m_+, m_-\}$ under some metric $g$. Suppose in addition that $(\Sigma^{15}, g)$ satisfies the $SC^p$-property at either the point $m_+$ or $m_-$. Then $\Sigma^{15}$ is diffeomorphic to $S^{15}$ or $\Sigma^{15}=[4064]\in bP_{16}=\mathbb{Z}_{8128}$.
\end{rem}

\begin{rem}
Very recently, Zhang and the second author \cite{TZ14} solved a
problem of B\'{e}rard-Bergery and Besse. That is, they showed that
every Eells-Kuiper quaternionic projective plane carries a
Riemannian metric with $SC^p$ property for certain point $p$. We
do not know how to improve our Proposition 3.2. For instance, we
do not know whether there is a metric on every Eells-Kuiper
quaternionic projective plane with not only the property in
Proposition 3.2, but also the $SC^p$-property.

\end{rem}

\begin{ack}
 The authors extend their heartfelt thanks to Professor Weiping Zhang for providing the reference \cite{DM90} and very useful comments on Moser's theorem. Thanks are also due to
 Dr. Jianqian Ge for his good ideas to prove Proposition \ref{propno} and Theorem 1.2.
\end{ack}

\end{document}